\documentclass{amsart}
\usepackage{amssymb,euscript,amsmath, mathrsfs}
\usepackage[dvips]{graphicx}
\usepackage[dvips]{color}

\newcounter{ENUM}
\newcommand{\itm}{\item}
\newenvironment{ilist}{\renewcommand{\theENUM}{\roman{ENUM}}\renewcommand{\itm}{\addtocounter{ENUM}{1}\item[(\theENUM)]}\begin{itemize}\setcounter{ENUM}{0}}{\end{itemize}}

\newenvironment{alist}[1][0]{\renewcommand{\theENUM}{\alph{ENUM}}\renewcommand{\itm}{\addtocounter{ENUM}{1}\item[\theENUM)]}\begin{itemize}\setcounter{ENUM}{#1}}{\end{itemize}}
\newenvironment{nlist}[1][0]{\renewcommand{\theENUM}{\arabic{ENUM}}\renewcommand{\itm}{\addtocounter{ENUM}{1}\item[(\theENUM)]}\begin{itemize}\setcounter{ENUM}{#1}}{\end{itemize}}

\newcommand{\margh}[1]{}

\input xy
\xyoption{all}
\CompileMatrices

\def\Z{{\mathbb Z}}

\def\A{{\mathbb A}}
\def\P{{\mathbb P}}

\def\fF{{\mathbb F}}
\def\F{{\mathscr F}}
\def\sG{{\mathscr G}}

\def\E{{\mathscr E}}

\def\sO{{\mathscr O}}
\def\I{{\mathscr I}}

\def\vp{\varphi}

\def\Hom{\operatorname{Hom}}

\def\Aut{\operatorname{Aut}}

\def\Tr{\operatorname{Tr}}
\def\Spec{\operatorname{Spec}}

\def\res{\operatorname{res}}

\def\rk{\operatorname{rk}}

\def\ord{\operatorname{ord}}

\numberwithin{equation}{section}
\newtheorem{thm}{Theorem}[section]
\newtheorem{prop}[thm]{Proposition}
\newtheorem{lem}[thm]{Lemma}
\newtheorem{cor}[thm]{Corollary}

\theoremstyle{definition}
\newtheorem{defn}[thm]{Definition}

\newtheorem{ex}[thm]{Example}
\newtheorem{sit}[thm]{Situation}
\theoremstyle{remark}
\newtheorem{notn}[thm]{Notation}
\newtheorem{rem}[thm]{Remark}
\newtheorem{warn}[thm]{Warning}

\begin{document}
\title{Logarithmic Connections With Vanishing $p$-Curvature}
\author{Brian Osserman}
\begin{abstract} We examine logarithmic connections with vanishing
$p$-curvature on smooth curves by studying their kernels, describing them in
terms of formal local decomposition. We then apply our 
results in the case of connections of rank $2$ on $\P^1$, classifying such 
connections
in terms of self-maps of $\P^1$ with prescribed ramification. 
\end{abstract}
\thanks{This paper was partially supported by fellowships from the National Science Foundation and Japan Society for the Promotion of Sciences.}
\maketitle

\section{Introduction}

We develop in this paper a basic theory of connections with simple poles 
and vanishing
$p$-curvature on smooth curves, and apply it to the case of rank $2$ vector
bundles on $\P^1$ to classify such connections completely in terms of
rational functions on $\P^1$ with prescribed ramification. 
Connections of this type are interesting in their own right, 
as demonstrated by a still-unsolved question of Grothendieck asking if a 
logarithmic connection on $\P^1$ in characteristic $0$ which has vanishing 
$p$-curvature when reduced mod $p$ for almost all primes $p$, must have
algebraic solutions; see \cite{ka2} for a discussion of the problem and 
solution for particular connections. However,
the immediate motivation for this paper was to use degeneration arguments
to conclude results on Frobenius-unstable vector bundles and the geometry of
the generalized Verschiebung on curves of genus $2$, as is carried out in
\cite{os11} and \cite{os9}, in much the same fashion as Mochizuki in
\cite{mo3}. 

We state our main theorem below in the case which will be of most use for
these applications, and which is simplest to state; however, its final
assertion
for connections can be obtained directly from Mochizuki's work. The most
general result, stated in Theorem \ref{con-main-nonconst}, may be used to 
conclude existence results for connections not treated by Theorem 
\ref{con-main}, nor by Mochizuki's results. Moreover, the classification of
our main theorem can also work in the other direction, obtaining
new results on self-maps of $\P^1$ via the classification here and work of
Mochizuki; for these last two applications, see \cite{os6}.

\begin{thm}\label{con-main} Fix an integer $n>0$, let $\delta = 0 \text{
or }1$ according to the parity of $n$, $d = \frac{n+\delta p}{2}-1$, 
and choose $P_1, \dots, P_n$ distinct points on $\P^1_k$, with $k$ an
algebraically closed field of characteristic $p>2$. Also fix $\E$ to be 
the vector bundle $\sO(\delta p -d) \oplus \sO(d)$. Then given an object
of the form $(\{\alpha_i\}_i, \bar{f})$, 
where the $\alpha_i$ are integers between $1$ and $\frac{p-1}{2}$, and 
$\bar{f}$ is a separable rational function on $\P^1$ of degree
$n (\frac{p-1}{2})+1-\sum_i \alpha_i$, and ramified to order at least
$p-2\alpha_i$ at each $P_i$, 
we can naturally construct a transport-equivalence class of 
connections on $\E$ with trivial determinant, vanishing $p$-curvature, 
simple poles at the $P_i$, and not inducing a connection on 
$\sO(d) \subset \E$. This association induces an injective map modulo the 
equivalence on rational functions of fractional linear transformation.

If further the $P_i$ are general, we obtain a bijective correspondence,
which holds even for first-order infinitesmal deformations.
In particular, the classes of such connections have no non-trivial 
deformations, and are counted by the recursive formula of 
\cite[Thm. 1.4]{os7}.
\end{thm}

We use throughout the standard terminological conventions for vector
bundles, connections, and $p$-curvature; see, e.g., \cite{ka1} for the
last two. Our methodology will be
to work primarily over an algebraically closed field, with periodic 
examinations of the generalization to first-order infinitesmal deformations. As such, we fix the following notation throughout.

\begin{notn} A `deformation' will always refer to a first-order infinitesmal deformation, and $\epsilon$ will always be a square-zero 
element. 
\end{notn}

We also specify the following terminology.

\begin{defn} We also say that $\nabla$ is a {\bf rational connection} on 
a smooth scheme if it is a connection on a dense open subset $U_{\nabla}$, 
but may have poles away from $U_{\nabla}$; we say that $\nabla$ is {\bf
logarithmic} if all such poles are simple.
\end{defn}

\begin{warn} We will refer to connections with {\bf trivial determinant} on
vector bundles $\E$ on $\P^1$ in the case that $p|\deg \E$, even if $\deg \E
\neq 0$, since in this case we have a unique canonical connection on 
$\det \E$, and can require that the determinant connection agree with it. This is a special case of the notion of a connection having $p$-trivial determinant, introduced in \cite{os6}.
\end{warn}

We begin in Section \ref{s-con-formal} with some calculations holding on any
smooth curve, the primary purpose of which is to show that a connection is
logarithmic with vanishing $p$-curvature if and only if everywhere formally
locally it decomposes as a direct sum of connections on line bundles. The
purpose of Section \ref{s-con-deform} is to re-establish the results of the
previous section for certain first-order infinitesmal deformations. Section
\ref{s-con-rank-2} develops simpler criteria in the special case of vector
bundles of rank 2, Section \ref{s-con-global} specializes further to the
case of vector bundles on $\P^1$, and Section \ref{s-con-maps} completes 
the classification in this situation in terms of self-maps of $\P^1$ with
prescribed ramification. 

The only similar work in the literature appears to be that
of Mochizuki, who proves a special case of the main results of this paper, 
in the situation of three poles on $\P^1$; in fact, he proves this result in
the more general context of $n$-connections over an arbitrary base, so our 
result (in the case of
three points) is simply the $n=0$ case of \cite[Thm. IV.2.3, p. 211]{mo3}.

The contents of this paper form a portion of the author's 2004 PhD thesis at
MIT, under the direction of Johan de Jong.

\section*{Acknowledgements}

I would like to thank Johan de Jong for his tireless and invaluable guidance. 

\section{Formal Local Calculations}\label{s-con-formal}

In this section, we make some basic observations about kernels of
connections with vanishing $p$-curvature and simple poles on smooth curves, 
and apply formal
local analysis to show that, formally locally, they may be split as a direct
sum of connections on line bundles; equivalently, they may be diagonalized
under transport.

We make the following definitions:

\begin{notn} We write $\E$ for a vector bundle of rank $r$ on $C$, and 
$\F$ for a vector bundle of the same rank on $C^{(p)}$. We also write 
$\vp$ for an injection $F^* \F \hookrightarrow \E$, and $\nabla$ for 
a connection on $\E$.
\end{notn}

\begin{defn} Given a vector bundle $\E$ of rank $r$ on $C$, we define a 
{\bf pre-kernel map} (to $\E$) to be a pair $(\F, \vp)$
with $\F$ locally free of rank $r$ on $C^{(p)}$, and $\vp:F^*\F \rightarrow
\E$ an injection. By abuse of terminology, we will refer to modification of 
$\vp$ by $F^*\Aut(\F)$ and $\Aut(\E)$ as {\bf transport}. 
\end{defn}

Note that a pre-kernel map induces a natural rational connection on $\E$ by
defining the sections $F^{-1} \F$ to be horizontal. 

\begin{defn} If a pre-kernel map $(\F, \vp)$ further has the property that
$\vp(F^{-1} \F)$ is the entire set of horizontal sections of the induced rational connection on $\E$, we say it is a {\bf kernel map}. 
\end{defn}

\begin{prop}\label{con-ker-lf} Let $C$ be a smooth curve over an 
algebraically closed field $k$, and $\E$
a vector bundle of rank $r$ on $C$. Then if we consider the operations of 
taking kernels of connections and of extending canonical connections of
Frobenius pullbacks, we deduce: 
\begin{ilist}
\itm There is a one-to-one correspondence between rational connections on 
$\E$ with vanishing $p$-curvature on one
side, and kernel maps $(\F, \vp)$ to $\E$ on the other, taken modulo automorphisms of $\F$.
\itm Under this equivalence, the poles of a connection are precisely the
points where $\vp$ fails to be surjective.
\itm Under this equivalence, transport of connections on $\E$ corresponds 
to changing $\vp$ by the corresponding automorphism of $\E$.
\end{ilist}
\end{prop}

\begin{proof} Let $\nabla$ be a rational connection on $\E$. Then since $C^{(p)}$ and $C$ are smooth curves, we find 
that $\E^\nabla$ and hence $F^* \E^\nabla$ are both vector bundles. 
Indeed, $F^* \E^\nabla$ is naturally
a subsheaf of $\E$, and can be understood concretely as the subsheaf
spanned by the kernel of $\nabla$ inside $\E$. We thus have a sequence
$0 \rightarrow F^* \E^\nabla \rightarrow \E \rightarrow \sG \rightarrow 0$
for some $\sG$ on $C$, and the inclusion map giving us the $\vp$ from 
statement (i). It follows from the Cartier isomorphism 
\cite[Thm. 5.1]{ka1} applied to the regular locus of $\nabla$ that 
$F^* \E ^ \nabla$ has rank $r$ if and only if $\sG$ is torsion if and only 
if $\nabla$ has vanishing $p$-curvature, and that in this case $\sG$ is 
supported at the poles of $\nabla$.
This yields one direction of (i), as well as (ii).

On the other hand, given a pre-kernel map $(\F, \vp)$, the induced
connection $\nabla$ satisfies $\E^\nabla \supseteq \vp(F^{-1}\F)$, and by the above, $\nabla$ has vanishing $p$-curvature. If we add the hypothesis that 
$\E^\nabla = \vp(F^{-1}\F)$, we find that $\nabla$ determines $\F$ uniquely, and 
$\vp$ up to automorphisms of $\F$ (note: not up to automorphisms of 
$F^* \F$, which will change $\nabla$), completing the 
proof of (i). Statement (iii) is now clear, completing the proof.
\end{proof}

We now carry out a straightforward calculation:

\begin{prop}\label{con-first-calc} Given $\E$ on $C$ and a pre-kernel map $(\F, \vp)$, let $\E$ and $F^{-1}\F$ be trivialized on an open subset $U$ of $C$, and suppose $\vp$ is given in terms of this trivialization by a matrix $S$. Then if $\nabla$ is the corresponding 
rational connection with vanishing $p$-curvature on $\E$, it has matrix 
$T=S(dS^{-1})$. Further, $\Tr(T) = -\frac{d \det(S)}{\det(S)}$.
\end{prop}

\begin{proof}This is straightforward linear algebra, using that $S$ is
generically invertible, the connection rule for $\nabla$, and that the image 
$\vp(F^{-1}\F)$ is in the kernel of $\nabla$ by definition.
\end{proof} 

We next move on to formal local analysis of the situation at points where the
determinant is not invertible (equivalently, points where the connection has
poles).

\begin{prop}Formally locally (that is, over $k[[t]]$), any $r \times r$ 
matrix of nonzero determinant:
\begin{ilist}
\itm may be put via left change of basis into the following form:
$$\begin{bmatrix}t^{e_1} & f_{12} & \cdots & & f_{1r} \\
0 & t^{e_2} & f_{23} & & \vdots \\
\vdots & & \ddots & \ddots & & \\
& & & & f_{(r-1)r} \\
0 & \cdots & & & t^{e_r}
\end{bmatrix}$$
where each $f_{ij}$ is a polynomial in $t$ of degree less than $e_j$;
\itm may, if one further allows right $p$th power change of basis,
be put into the above form, with the further requirement that the 
$f_{ij}$ do not have any terms with exponent congruent to $e_i$ modulo $p$.
\end{ilist}
\end{prop}

\begin{proof} The form of (i) may be obtained by standard row reduction techniques. For (ii), we remove the terms congruent to $e_i$ mod $p$ from each $f_{ij}$ using $p$th-power column reduction.
\end{proof}

\begin{rem}Note that unlike form (i) of the preceding proposition, form 
(ii) is not unique. In particular, conjugation by permutation matrices is 
always allowed, and could be used to rearrange the coefficients of a
diagonal matrix; this could not be accomplished using row
reduction alone.
\end{rem}

\begin{prop}\label{con-formal} For a pre-kernel map $\vp$ given on some open subset by 
$S=(a_{ij})$, the following are equivalent:
\begin{alist}
\itm $\vp$ corresponds to a logarithmic connection with 
vanishing $p$-curvature;
\itm formally locally everywhere (equivalently, everywhere where the 
map fails to be invertible), $S$ is transport-diagonalizable, with all
diagonal coefficients having order of vanishing strictly less than $p$;
\itm formally locally everywhere (equivalently, everywhere where the 
map fails to be invertible), when $S$ is placed in the form of the preceding 
proposition, all $f_{ij}=0$ and all $e_i$ are strictly less than $p$.
\end{alist}
\end{prop}

\begin{proof} First note that the condition that $\vp$ correspond to a
$\nabla$ with vanishing $p$-curvature and at most simple poles is clearly
transport-invariant. We do the difficult direction first; namely, showing
that a) implies c). For notational 
convenience, we prove this inductively on the
rank $r$. The base case is $r=1$, where the connection corresponding to
$a_{11}$ is simply $-\frac{da_{11}}{a_{11}}$, which always has at most
simple poles. The condition that $e_1 < p$ comes from the fact that if $e_1
\geq p$, and we have $S=\begin{bmatrix}t^{e_1}\end{bmatrix}$, and
$T=\begin{bmatrix}-e_1 t^{-1}dt\end{bmatrix}$, then $t^{e_1-p}$ will also be
a horizontal section formally locally, but is not in the image of $\vp$.
Here we are using that in characteristic $p$, because a connection on $C$ is
$\sO_{C^{(p)}}$-linear, formation of the kernel of a
connection commutes with completion.

For the induction step, we first transport $S$ formally locally into 
the form described in part (ii) of the previous proposition; this is in 
particular
upper triangular, and noting that once $S$ is upper triangular $T$ is also
upper triangular, we can (formally locally) restrict to the first $r-1$ 
rows and columns of $T$ to get a connection with vanishing $p$-curvature and simple poles in
rank $r-1$, which is clearly already in the form of the previous
proposition. Thus, by the induction hypothesis our entire $r \times
r$ matrix will look like:

$$\begin{bmatrix}t^{e_1} & 0 & \dots & 0 & f_1 \\
0 & \ddots & \ddots & \vdots & \vdots \\
\vdots & \ddots & t^{e_{r-2}} & 0 & f_{r-2}\\
\vdots & & \ddots & t^{e_{r-1}} & f_{r-1}\\
0 & \dots & \dots & 0 & t^{e_r} \end{bmatrix}$$

We wish to show that the $f_i$ are all $0$, and $e_r < p$. Computing the
associated connection matrix, we see that we need only consider the last column, which will
have simple poles only if for each $i<r$, the function $e_i f_i - f_i' t$ vanishes to order at
least $e_r$; since $f_i$ has degree less than $e_r$ by hypothesis, this
difference must be $0$. But it is clear that terms will cancel in a given
degree if and only if the degree is congruent to $e_i$ mod $p$, and also by
hypothesis each $f_i$ has no terms in degree congruent to $e_i$ mod $p$. We
conclude that each $f_i=0$, as desired. Lastly, the condition that $e_r < p$
follows from the necessity of the image of $\vp$ to contain the kernel of 
$\nabla$ just as it did in the rank $1$ base case.

Now, c) implies b) is trivial, so we just need to show that b) implies
a). If $S$ is formally locally diagonalizable, as long as the $e_i$ are 
less than $p$ the diagonalized map corresponds to a connection with simple 
poles and vanishing
$p$-curvature, and since this is a transport-invariant property, $S$ must
have as well.
\end{proof}

Because under these equivalent conditions, all $e_i < p$, we note that in fact it 
is only necessary to use constant column operations in our formal
local transport-diagonalization procedure, so we conclude:

\begin{cor}\label{con-ct} A pre-kernel map $\vp$ given on some open subset 
by $S=(a_{ij})$
corresponds to a logarithmic connection with vanishing $p$-curvature
if and only if at each point where $\vp$ fails to be surjective, for
$t$ a local coordinate at that point, there exist constants $c_{ij}$ for
all $0<i<j \leq r$ and a formal local invertible $M$ such that $MSU(c_{ij})$ 
is diagonal with $t^{e_i}$ as its diagonal coefficients, and all $e_i<p$, 
where $U(c_{ij})$ is the upper triangular matrix having $1$'s on 
the diagonal and given by the $c_{ij}$ above the diagonal.
\end{cor}

We may also phrase this last result purely in terms of connections:

\begin{cor}\label{con-diag} A rational connection (having at least one pole) 
is logarithmic
with vanishing $p$-curvature if and only if, formally locally at every 
pole, the connection
may be transported so as to have diagonal matrix with each diagonal entry of
the form $e_i t^{-1} dt$, with $e_i \in \fF_p$, or equivalently, if the
connection decomposes formally locally as a direct sum of connections on
line bundles, each with logarithmic poles and vanishing $p$-curvature.
\end{cor}

\begin{proof} The only if direction follows immediately from our prior work:
by Corollary \ref{con-ct}, the kernel map is formally locally diagonalizable 
with diagonal entries $t^{e_i}$, and by Proposition \ref{con-first-calc} we
see that this gives a connection of the desired form. Conversely, one
computes directly that given a diagonal connection as described, the kernel
mapping may be given explicitly by a diagonal matrix with $t^{e_i}$ on the
diagonals (where $0 \leq e_i <p$), and is in particular of full rank,
implying that the $p$-curvature of the connection vanishes.
\end{proof}

As another corollary, we can put together the preceding 
propositions to get the following relationship between $\det(S)$, $\Tr(T)$, 
the $e_i$, and the eigenvalues of the residue matrix $\res_t T$: 

\begin{cor}\label{con-det-formula} With the notation of Proposition
\ref{con-first-calc}, if $\nabla$ is logarithmic with vanishing $p$-curvature, $\res_t T$ is
diagonalizable (in the usual sense), with eigenvalues given
as the $e_i$ mod $p$. The determinant satisfies $\ord_t \det(S) = \sum _i e_i 
\equiv \Tr(\res_t T) \pmod{p}$, but moreover, if we have the $e_i$ only in 
terms of 
their reductions $\bar{e}_i$ mod $p$, we also have the formula $e_i =
<\bar{e}_i>$, and hence $\ord_t \det(S) 
= \sum _i <\bar{e}_i>$, where $<a>$ for any $a \in \Z/p\Z$ denotes the unique
integer representative for $a$ between $0$ and $p-1$. Finally, transport of
$T$ along an automorphism conjugates $\res_t T$ by the same automorphism.
\end{cor}

From here on we assume that we are in the following situation:

\begin{sit}\label{con-sit} Our connection $\nabla$ is logarithmic, with 
vanishing $p$-curvature. At every pole of $\nabla$, we suppose that the
$e_i$ of Corollary \ref{con-diag} are all non-zero.
\end{sit}

The non-vanishing conditions on the $e_i$ will come into play only when we
attempt to study deformations of connections. 

\begin{rem} Although the determinant of a connection, and in particular its triviality, is well-determined under transport
equivalence only globally on a proper curve, the trace of the residue of a connection is always invariant. Indeed, 
an automorphism given locally by a matrix $S$ will act on a connection matrix
$T$ by $T \rightarrow S^{-1} T S + S^{-1} dS$, and invariance of the trace
of the residue follows from the invertibility of $S$.
\end{rem}

\begin{rem}We cannot expect such nice behavior when we weaken the hypothesis
that $\nabla$ have only simple poles. First of all, it is easy enough to
write down examples of higher order poles, as soon as the rank is higher 
than one. In this
situation, the relationship between the order of the determinant and the
order of the poles is much less clear-cut. Moreover, it is easy to check
that for rank higher than one, at a point with poles of order greater than
one, the residue itself is no longer well-defined under transport.
\end{rem}

\section{Generalization to $k[\epsilon]$}\label{s-con-deform}

The aim of this section is to generalize the results of the previous section
to the case where we have changed base to $\Spec k[\epsilon]$. It turns out
that the most difficult part of this is to show that the kernel of an 
appropriate deformation of a connection as in Situation \ref{con-sit} will 
give a deformation of the kernel of the original connection. 
We proceed in several steps. We first pin down our situation and notation:

\begin{sit}\label{con-deform} We suppose that $C$ is obtained from a smooth
proper curve $C_0$ over $k$ via change of base to $\Spec k[\epsilon]$, and
similarly for a vector bundle $\E$ on $C$ from $\E_0$ on $C_0$. We have a 
connection $\nabla$ on $\E$, and $\nabla_0$ is the induced connection on 
$\E_0$.
\end{sit}

\begin{notn} If $D$ is the divisor of poles of $\nabla_0$, so that
$\nabla_0$ takes values in $\E_0 \otimes \Omega^1_{C_0}(D)$, then we denote by
$\I_{\nabla_0}$ the $\sO_{C_0}$-submodule of $\E_0 \otimes \Omega^1_{C_0}(D)$ 
generated by the image of $\nabla_0$.
\end{notn}

Our first goal will be to show that in our situation, with very minor
additional hypotheses, $\E^{\nabla}$ is a deformation of $\E_0^{\nabla_0}$. Specifically:

\begin{prop}\label{con-ker-def} Suppose $\nabla$ is a logarithmic connection
with vanishing $p$-curvature, with $\nabla_0$ having poles wherever $\nabla$
does, and such that all the $e_i$ of Corollary \ref{con-diag} applied to 
$\nabla_0$ are non-zero. Then:
\begin{ilist}
\itm $\E^{\nabla}$ is locally free on $C^{(p)}$ of rank equal to $\rk \E$; 
\itm the natural map $\E^{\nabla}/\epsilon \E^{\nabla}
\rightarrow (\E_0)^{\nabla_0}$ is an isomorphism.
\end{ilist}
\end{prop}

\begin{proof} We first claim that to prove (i), it will suffice to show that
$\E^{\nabla}/\epsilon \E^{\nabla}$ is torsion-free over $C_0$. Indeed, one can check that in our situation of $k[\epsilon]/(\epsilon^2)$, it is enough see that $\E^{\nabla}/\epsilon \E^{\nabla}$
and $\epsilon \E^{\nabla}$ are both locally free over $C_0$, of rank equal
to the rank of $\E$, without any {\it a priori} hypotheses on the natural map between them. Now, on the open subset of $C$ on which $\nabla$ is 
regular, it follows easily from the Cartier isomorphism that
we have that $\E^{\nabla}$ is locally free of the
correct rank, and so then are $\E^{\nabla}/\epsilon \E^{\nabla}$ and
$\epsilon \E^{\nabla}$. The required rank condition will thus follow
automatically if we can show that both these sheaves are locally free on all
of $C_0$, which is a smooth curve; this reduces the problem to showing that 
both these sheaves are torsion-free. Finally, $\epsilon \E^{\nabla}$ is a
subsheaf of the locally free sheaf $\E$ and hence torsion-free, so we obtain
the desired reduction of (i) to showing that $\E^{\nabla}/\epsilon
\E^{\nabla}$ is torsion-free over $C_0$.

We now reduce both (i) and (ii) down to a certain divisibility lemma. Both
statements are local on $C$, so we make our analysis entirely on stalks,
letting $P$ be an arbitrary point of $C$. Locally, $\E$ is free, 
so we can pick a splitting map $\E/ \epsilon \E \rightarrow \E$. 
We can then write $\nabla = \nabla_0 +
\epsilon \nabla_1$, and it makes sense to view both $\nabla_0$ and
$\nabla_1$ as taking values in $\E /\epsilon \E$ (since this is naturally
isomorphic to $\epsilon \E$). The basic observation is that $\nabla_1$ must
take values in $\I_{\nabla_0}$: indeed, it may have simple poles only where
$\nabla_0$ does, so it takes values in $\E_0 \otimes \Omega^1_{C_0}(D)$, and 
by Corollary \ref{con-diag}, we see by the hypothesis 
that all the $e_i$ are non-zero that $\I_{\nabla_0}$ is all of 
$\E_0 \otimes \Omega^1_{C_0} (D)$.

We first consider (i): since we are checking that $\E^{\nabla}/\epsilon
\E^{\nabla}$ has no torsion as a module over $\sO_{C^{(p)}}$, we need only
consider multiplication by $f \in \sO_{C,P}$ such that $df=0$. We must show
that given $s \in \E^{\nabla}_P$ with $fs \in \epsilon \E^{\nabla}_P$, we
must have $s \in \epsilon \E^{\nabla}_P$. If we write $fs = \epsilon s'$, with $s'
\in \E^{\nabla}_P$, and $s' = s'_1 + \epsilon s'_2$, then we see that $f |
s'_1 \in \E_{0,P}$, and it will suffice to show we can choose $s'_2$ so that
$f| s'_2$ as well, since then we can divide through by $f$ to write $s$ as 
$\epsilon$ times an element of $\E^{\nabla}_P$. Since the value of $s'$ is
only relevant modulo $\epsilon$, we may replace $s'_2$ by any element which
keeps $s'$ in the kernel of $\nabla$. Now, we have
$0 = \nabla(s') = \nabla_0 (s'_1)+\epsilon (\nabla_1(s'_1) +
\nabla_0(s'_2))$, and since $df=0$, it follows that $f| \nabla_1(s'_1)$. 
Because both $\nabla_1$ and $\nabla_0$ take values in $\I_{\nabla_0}$, we must
have $f| \nabla_0(s'_2)$ in $\I_{\nabla_0}$, and the divisibility lemma which
follows completes the proof, taking $s'_2$ as our $s$ in the lemma, and 
obtaining our new $s'_2$ as the lemma's $fs'$.

Next, we wish to reduce (ii) down to the same lemma. Having already
completed (i), we may assume that $\E^{\nabla}$ is locally free, with rank
equal to $\rk \E$. It follows that $\E^{\nabla}/\epsilon \E^{\nabla}$ is
locally free of the same rank on $C_0$, as is $(\E_0)^{\nabla_0}$ by 
Proposition \ref{con-ker-lf}. It therefore suffices to show that the natural 
map is a
surjection in order to conclude that it is an isomorphism. Let $s_0$ be a
section of $(\E_{0,P})^{\nabla_0}$; we need only lift it to a section $s \in
\E^{\nabla}_P$. Moreover, we know that we can do so generically, since we
have the Cartier isomorphism away from the poles of $\nabla$ by 
\cite[Thm. 5.1]{ka1}. Therefore, there exists some $f$ such that $fs_0$ 
lifts to a section $s$ of $\E^{\nabla}_P$; as before, we are working over
$\sO_{C^{(p)}}$, so as an element of $\sO_{C,P}$, we have $df=0$. But now we 
find
ourselves in the same situation as before: if $s = s_1 + \epsilon s_2$, we
have that $f|s_1$, we want $f$ to divide $s_2$, and we may modify $s_2$
arbitrarily as long as $s$ remains in the kernel of $\nabla$. Thus by the
same argument as for (i), we reduce to our divisibility lemma.
\end{proof}

\begin{lem}\label{con-pdivide} We continue with the hypotheses of the 
previous proposition. Given $f \in \sO_{C,P}$ for some $P \in C$, with 
$df=0$, and $s$ in the stalk $\E_{0,P}$ with $f | \nabla_0(s)$ in the 
stalk $\I_{\nabla_0, P}$, then there exists $s' \in \E_{0,P}$ with 
$\nabla_0(fs') = f\nabla_0(s')=\nabla_0(s)$.
\end{lem}

\begin{proof} Under our hypotheses on $\nabla_0$, which allow us to invoke
Corollary \ref{con-diag}, the proof is straightforward. We first
prove the result formally locally. In this setting, we claim it is enough to
handle the case $f = t^p$: in general, write $f = (t^p)^i u$ for some $i
\geq 0$ and some unit $u$; certainly, if we have handled the case of $t^p$, we
can inductively ``divide out'' $i$ times by $t^p$, and then since $u$ is a
unit, we can simply set $s' = u^{-1} s$. But for $f=t^p$, we simply carry
out a direct computation; the diagonalizability obtained from Corollary 
\ref{con-diag} expresses the connection formally locally as a direct sum of
connections on line bundles, so it suffices to work with rank one, and a
connection of the form $\nabla_0(s) = ds + e t^{-1} dt$, for some $e \in
\fF_p$; our $\I_{\nabla_0}$ in this context is simply everything of the form
$\sum_{i \geq -1} a_i t^i dt$. If we write $s = \sum _{i \geq 0} a_i t^i$, we 
get $\nabla_0(s) = \sum _{i \geq 0} (i + e) a_i t^{i-1} dt$; this is divisible
by $t^p$ in $\I_{\nabla_0}$ if and only if $(i+e) a_i = 0$ for all $i <p$.
Now, for any $i<p$ with $i+e=0$, we can replace $a_i$ with $0$ without
changing $\nabla_0(s)$, and for all other $i$, we must have $a_i=0$ to start
with. Hence, we see that we can modify $s$ in degree $p-e$, if necessary, so
that all $a_i = 0$ for $i<p$, and we can then obtain our $s'$ as $t^{-p}$
times our modified $s$. 

This gives the formal local result, but it is now easy enough to conclude
the desired Zariski-local statement. We have $s-fs'$ in the kernel of
$\nabla$, and because in characteristic $p$ formation of the kernel of a
connection commutes with completion, we
can write $s-fs' = \sum _i f_i s_i$ where $s_i \in \E_{0,P}^{\nabla}$ and
$f_i \in k[[t]]$. But by definition, we can approximate the $f_i$ to
arbitrary powers of $t$ by elements of $\sO_{C,P}$; if we let $f'_i$
approximate the $f_i$ to order at least $\ord_t f$, we find that $f$ must
divide $s - \sum_i f'_i s_i$, so we can set our desired Zariski-local section 
to be $\frac{1}{f}(s-\sum_i f'_i s_i)$.
\end{proof}

We now know the correct conditions for connections over $k[\epsilon]$.
Specifically, after this section, whenever we are over $k[\epsilon]$, we 
assume we have:

\begin{sit}\label{con-sit-2} Our connection $\nabla$ is logarithmic, with 
vanishing $p$-curvature. If $\nabla_0$ is the connection obtained modulo
$\epsilon$, then every pole of $\nabla$ must also be a pole of $\nabla_0$,
and we suppose that the $e_i$ of Corollary \ref{con-diag} as applied to
$\nabla_0$ are all non-zero.
\end{sit}

Finally, we are ready to conclude:

\begin{cor}\label{con-ct2} Corollary \ref{con-ct} holds even over 
$k[\epsilon]$; more precisely, a pre-kernel map $\vp$ as in Proposition
\ref{con-ker-lf}, given by $S=(a_{ij})$ on some open subset which contains
every point where $\vp$ fails to be surjective,
corresponds to a connection satisfying the conditions of
Situation \ref{con-sit-2}
if and only if at each point where $\vp$ fails to be surjective, for
$t$ a local coordinate at that point, there exist constants $c_{ij}\in
k[\epsilon]$ for
all $0<i<j \leq r$ and a formal local invertible $M$ such that $MSU(c_{ij})$ 
is diagonal with $t^{e_i}$ as its diagonal coefficients, and all $e_i<p$, 
where $U(c_{ij})$ is the upper triangular matrix having $1$'s on 
the diagonal and given by the $c_{ij}$ above the diagonal.
\end{cor}

\begin{proof} We first note that given an $S$, the calculation of Proposition 
\ref{con-first-calc} is still valid because $S$ and hence $\det S$ is still 
generically invertible. Hence, as before it is clear that if the desired 
$M, U(c_{ij})$
exist, then $S$ corresponds to a connection of the desired type. Conversely,
given such a connection, since $S$ describes the kernel of our connection,
by the previous proposition, we find that we have
an $S$ which agrees modulo $\epsilon$ with the $S_0$ obtained from taking the
connection modulo $\epsilon$; we can then apply Corollary \ref{con-ct}
to conclude that formally locally on $C$ there is an invertible $M_0$ and a
$U_0(\bar{c}_{ij})$, both over $k$, such that $S' = M_0 S U_0(\bar{c}_{ij})$
is of the desired form modulo $\epsilon$. Thus, we can write 
$$ S' = \begin{bmatrix} t^{e_1}+\epsilon f_{11} & \dots & \epsilon f_{1r} \\
\vdots & \ddots & \vdots \\
\epsilon f_{r1} & \dots & t^{e_r} + \epsilon f_{rr} \end{bmatrix}$$

One then checks that $T' = S'(dS'^{-1})$ is given by 
$$\begin{bmatrix} -e_1 t^{-1} + \epsilon t^{-e_1-1} (e_1 f_{11} - 
t f_{11}') & \dots & \epsilon t^{-e_r-1} (e_1 f_{1r} 
- t f_{1r}') \\
\vdots & \ddots & \vdots \\
\epsilon t^{-e_1-1} (e_r f_{r1} - t f_{r1}') 
& \dots & -e_r t^{-1} + \epsilon t^{-e_r-1} (e_r f_{rr} - 
t f_{rr}') \end{bmatrix} dt$$

Thus, in order to have simple poles, it is necessary and sufficient that
$\ord_t (e_i f_{ij} - t f_{ij}') \geq e_j$ for all $i,j$. But this is
precisely the condition required to be able to remove all the $f_{ij}$ via
row and (constant) column operations, since the inequality above implies
that all terms of $f_{ij}$ in degree $\ell$ must vanish for $\ell < e_j$,
unless $\ell = e_i$. Constant column operation can remove the terms of
degree $e_i$ from each $f_{ij}$, and then we have that $\ord_t f_{ij} \geq
e_j$, so row operations can remove the $f_{ij}$, as desired.
\end{proof}

\section{Applications to Rank $2$}\label{s-con-rank-2}

As our case of primary interest, we will develop the theory further in the
case of vector bundles of rank $2$ and connections $\nabla$ whose residue at 
all poles has trace zero. Note that in this case, at any pole the $e_i$ of
Corollary \ref{con-diag} satisfy $e_1 = - e_2$, and in particular are
automatically both non-zero as required in Situation \ref{con-sit}. 
We will work simultaneously over $k[\epsilon]$, assuming in this case the
conditions of Situation \ref{con-sit-2}. In this scenario, we define:

\begin{defn} Given $f \in A[[t]]$, we say that $\ord_t f = e$ if and only if
the first non-zero coefficient of $f$ is the coefficient of $t^e$. If 
further this first
non-zero coefficient is a unit in $A$, we say that $f$ {\bf vanishes
uniformly to order $e$} at $t=0$.
\end{defn}

Now, the kernel map $\vp$ associated to any connection $\nabla$ is given 
locally by a matrix 
$S=\begin{bmatrix}g_{11} & g_{12} \\
g_{21} & g_{22}\end{bmatrix}$. The corresponding connection $\nabla$ is 
then given locally
by a matrix $T$, which Proposition \ref{con-first-calc} allows us to 
write explicitly as

\begin{equation}\label{con-2x2T}
T=\frac{1}{\det S}\begin{bmatrix}(dg_{12}) g_{21} - (dg_{11}) g_{22} &
(dg_{11}) g_{12} - (dg_{12})g_{11} \\
(dg_{22})g_{21}-(dg_{21})g_{22} & 
(dg_{21})g_{12}-(dg_{22})g_{11} \end{bmatrix}
\end{equation}

Corollary \ref{con-det-formula} tells us that 
the simple poles of the connection will occur at precisely the places where 
$\det S$ vanishes, and that this will always occur to order precisely $p$.
Over $k[\epsilon]$, Corollary \ref{con-ct2} implies that the determinant 
will vanish uniformly to order $p$.
As before, choose a point where this is the case, and let $t$ be a local
coordinate at that point. Denote by $e_{ij}$ the order at $t$ of $g_{ij}$.
We will develop more precisely the criterion for $S$ to correspond to $T$
(that is, for the image of $S$ to contain the kernel of $T$), and for $T$ to
have simple poles. We find:

\begin{prop}\label{con-crit2} Over $k$ (respectively, $k[\epsilon])$, assuming 
that $\det S$ vanishes (uniformly) to order $p$ at $t=0$, for $S$ 
to correspond to a connection $T$ with a simple pole at
$t=0$ and vanishing $p$-curvature, it is necessary and sufficient that there
exists a $c_t$ such that after $S$ is replaced by 
$S'= S \begin{bmatrix}1 & -c_t \\ 0 & 1\end{bmatrix}$, we have:
$$\min\{\ord_t g_{11}, \ord_t g_{21}\} + \min \{\ord_t g_{12}, \ord_t
g_{22}\} \geq p.$$
Over $k$, this may be stated equivalently, after $S$ is replaced by
$S'$, as the condition that order of vanishing at $t=0$ be greater than or 
equal to $p$ for
all of $g_{11} g_{22}, g_{21} g_{12}, g_{11} g_{12}, g_{21} g_{22}$.
\end{prop}

\begin{proof}
First, if $S$ corresponds to a connection with a simple pole at $t=0$,
by Corollary \ref{con-ct} (respectively, Corollary \ref{con-ct2}) we have a 
$c_{12}$ such that 
$MS \begin{bmatrix}1 & c_{12} \\ 0 & 1\end{bmatrix}$ is diagonal with powers
of $t$ on the diagonal, and $M$ is a formal local invertible matrix. Letting
$c_t=-c_{12}$, we replace $S$ by $S'$, and are simply saying that $MS$ is
diagonal with powers of $t$ on the diagonal, say $t^e$ and $t^{p-e}$. Multiplying by $M^{-1}$, we trivially obtain the desired conditions on
the $g_{ij}$.

Conversely, suppose that the required $c_t$ exists, and we have replaced $S$
by $S'$. Let $e = \min\{\ord_t g_{11}, \ord_t g_{21}\}$. By hypothesis,
$\min\{\ord_t g_{12}, \ord_t g_{22}\} \geq p-e$. Thus,
we can write $S$ as $(m_{ij})D(t^e, t^{p-e})$ for some $m_{ij}$ regular at
$t=0$, and once again the condition that $S$ has determinant vanishing
uniformly to order $p$ at $t=0$ implies that $\det (m_{ij})$ is a unit, and
hence that $(m_{ij})$ is invertible and may
be moved to the other side, letting us apply Corollary \ref{con-ct}
(respectively, Corollary \ref{con-ct2}) to conclude
that $S$ corresponds to a connection with a simple pole at $t=0$ and
vanishing $p$-curvature, as desired.
\end{proof}

\begin{rem}This criterion looks rather asymmetric on the face of it, but
note that locally one may always conjugate by $\begin{bmatrix}0 & 1 \\1 & 0
\end{bmatrix}$ to switch the rows and columns, after which application of
the above criterion gives equivalent conditions in terms of
subtracting the right column from the left rather than vice versa. We will 
refer to this as the {\bf mirror criterion}.
\end{rem}

\begin{rem} Our initial description of $c_t$ was that there exist an 
invertible $M$ such that $M S \begin{bmatrix}1 & -c_t \\ 0 & 1\end{bmatrix}$ 
is diagonal, from which it immediately follows that $c_t$ is independent of
transport of $S$ via left multiplication. However, if we multiply by 
some invertible column-operation matrix $N$ on the 
right, we will need to determine how to ``move'' this action over to the left,
which is in general not a simple matter, and can result in substantial
changes to the behavior of the $c_t$. This is a rather ironic situation,
since it is the right multiplication which leaves the corresponding
connection unchanged, and the left multiplication which applies automorphism 
transport to it. In any case, we will at least be able to characterize 
exactly how the $c_t$ can change under global right multiplication in most
cases on $\P^1$. 
\end{rem}

\section{Global Computations on $\P^1$}\label{s-con-global}

Throughout this section, let $\E$ be $\sO(\delta p -d) \oplus \sO(d)$ on 
$\P^1$, where $\delta = 0 \text{ or } 1$, and $\delta p < 2d$.
We set up the basic situation to be used in this section and the
next, and then classify in Proposition \ref{con-constant} an ``easy case''
for the connections we wish to study, which will not be relevant to Theorem 
\ref{con-main}, but which we include nonetheless for the sake of 
completeness. Let $t_1, \dots, t_n$ be local coordinates at $n$ distinct 
points on $\P^1$; without loss of generality, write $t_i = x - \lambda_i$, 
where $x$ is a coordinate for some $\A^1 \subset \P^1$ containing the 
relevant points,
and let $c_i$ be the $c_t$ of Proposition \ref{con-crit2} for each $t_i$.
If $\nabla$ is a connection on $\E$ with vanishing $p$-curvature and simple 
poles at the $\lambda_i$, $F^*\E^{\nabla}$ must
have degree $-p(n - \delta)$, so it must be of the form 
$\sO(-mp) \oplus \sO((m-n+\delta) p)$
for some integer $m$ (without loss of generality, say $m \leq n-\delta-m$), 
and because it must map with full rank to $\sO(\delta p -d) \oplus \sO(d)$, 
we find that
we also must have $(m-n+\delta)p \leq \delta p-d$, $-mp \leq d$, which 
gives us 
$-d \leq mp \leq np-d$. We now fix some choice of $m$, and consider
possibilities for the kernel map $\vp$ corresponding to such a $\nabla$.

We may write $\Hom(F^* \E^{\nabla}, \E)$ as
$$\begin{bmatrix}\sO((m+\delta)p-d) & \sO((n-m)p-d) \\ \sO(mp+d) & 
\sO((n-\delta-m)p+d) \end{bmatrix}$$
The matrix $S$ can therefore be written with coefficients $g_{ij}$ being 
polynomials
in $x$ of the appropriate degrees, with products along both the diagonal and
antidiagonal having
degree bounded by $np$. Moreover, there are $n$ points where the determinant
must vanish to order $p$, so up to scalar multiplication, the determinant
must be $\prod _i (x-\lambda_i)^p$. Global transport of our kernel map
corresponds to left multiplication by matrices of the form
$$\begin{bmatrix}\sO(0) & \sO(2d-\delta p) \\ \sO(\delta p -2d) & \sO(0)
\end{bmatrix}$$
and right multiplication by
$$F^* \begin{bmatrix}\sO(0) & \sO(n-\delta-2m) \\ \sO(2m-n+\delta) & \sO(0)
\end{bmatrix}$$

Then we have: 

\begin{prop}\label{con-cimove}Although the $c_i$ are not invariants of a connection, for the
most part they change predictably under transport of their kernel maps. It is 
always possible to
scale them simultaneously. It is also possible to translate them
simultaneously by
(any constant times) $\lambda_i^{pj}$ for any $j$ between $0$ and
$n-\delta-2m$. If $m<n-\delta-m$, the $i$ for which the $c_i$ are 
uniquely defined do not change 
under transport, and the above modifications are the only possible ones for 
these $c_i$.
\end{prop}

\begin{proof}We make use only of the criterion of Corollary \ref{con-ct}
(recalling that the $c_i$ were by definition the negatives of the constants
arising there).
We first show that the asserted modifications are possible. If we begin with
$S$, and at each $\lambda_i$ an $M_i$ and upper triangular $U(-c_i)$ with 
$M_i S U(-c_i)$ diagonal, we can transport $S$ to simultaneously scale the 
$c_i$ by any $\mu$ simply by replacing $S$ by 
$S D(1,\mu)=S \begin{bmatrix}1 & 0 \\ 0 & \mu \end{bmatrix}$, 
$U(-c_i)$ by $D(1,\mu^{-1})U(-c_i)D(1,\mu) = U(-\mu c_i)$, and $M_i$ by
$D(1,\mu^{-1}) M_i$, whereupon our original diagonal matrix is conjugated 
by $D(1,\mu)$.

Next, translation of all $c_i$ by $\mu \lambda_i^{pj}$
is accomplished simply by right multiplication of $S$ by 
$U(\mu x^{pj})$: at each $\lambda _i$, we can write $x^{pj} =
\lambda_i ^{pj} + (x^j-\lambda_i^j)^p$, and then if 
$M_i S U(-c_i)=D(d_1, d_2)$ was diagonal, it follows that 
$M_i (S U(\mu x^{pj})) U(-c_i-\mu \lambda_i^{pj}) = M_i S U(-c_i) 
U(\mu (x^j-\lambda_i^j)^p) =
\begin{bmatrix}d_1 & \mu d_1 (x^j-\lambda_i^j)^p \\ 0 & d_2\end{bmatrix}$. 
Now, since $\ord_{\lambda_i} d_2 < p$, we can multiply $M_i$ on the left by
$U(-\mu \frac{d_1 (x^j-\lambda_i^j)^p}{d_2})$ to recover the initial diagonal
matrix, so we see that $c_i + \mu \lambda_i^{pj}$ has taken the role of
$c_i$, as desired.

Lastly, when $m<n-\delta -m$, we simply need to verify that the above 
cases are the
only possible forms of transport that can affect the $c_i$: we have seen
that only right multiplication can affect the $c_i$, and when 
$m<n-\delta-m$, the
only matrices we can right multiply by are upper triangular with scalars on 
the diagonal and inseparable polynomials of degree $\leq (n-\delta-2m)p$ 
in the 
upper right. These are generated by the two cases above together with
$D(\mu, 1)$, but $D(\mu, 1)=D(\mu, \mu) D(1, \mu^{-1})$, and the $D(\mu,
\mu)$ can be commuted to the left and absorbed into $M$. In particular, 
all methods 
of acting on the $c_i$ change them invertibly, so whether or not they are
uniquely determined is tranport-invariant as long as $m<n-\delta-m$.
\end{proof}

\begin{ex} When $m=n-\delta-m$, it is not true that the $c_i$ behave well 
under
transport, and they may even go from uniquely determined to arbitrary and
back. For instance, consider a diagonal matrix vanishing to order $e<p/2$,
$p-e$ along the diagonal at a chosen point $\lambda_i$. In this case, $c_i$
is well-determined as $0$, since if we multiplied by any $U(c_i)$ with $c_i
\neq 0$, we would have that the product of the entries on the top row of our
matrix only vanished to order $2e < p$. But because $m=n-\delta-m$, we can
right-multiply by $\begin{bmatrix}0 & 1 \\ 1 & 0 \end{bmatrix}$ to switch
the columns of our matrix, at which point $c_i$ may be chosen arbitrarily,
because $2p-2e > p$.
\end{ex}

Before moving on to the next results, we fix some combinatorial notation
which will come up as soon as we attempt to count classes of
connections. 

\begin{notn}\label{con-np} For a given $p$, $n$, and $s$, denote by $N_p(n,s)$ the number
of monomials of degree $s$ in $n$ variables subject to the restriction
that each variable occur with positive exponent strictly less than $p$. Also
denote by $N_p^D(n,s)$ the number of such monomials in which exactly $D$
variables appear with degree less than $p/2$.
\end{notn}

We give explicit formulas for these numbers:

\begin{lem}We have:
$$N_p(n,s)=\sum_{i=0}^n (-1)^i \binom{n}{i}\binom{s-i(p-1)-1}{n-1}.$$ In
particular, for any fixed $n$ this is expressed
by the $j$th of $n+1$ polynomials, each of combined degree $n-1$ in $s$ and 
$p$, 
with $j$ being the largest integer ($\leq n+1$) such that $j(p-1) \leq s-n$.
We also have 
$$N^D_p(n,s) = \binom{n}{D}\sum_{j=0}^s N_{(p+1)/2}(D,j)N_{(p+1)/2}(n-D, s -
j - (n-D)\frac{p-1}{2}).$$
\end{lem}

\begin{proof} The formula for $N_p(n,s)$ is obtained by the
inclusion-exclusion principle, looking at which variables occur with 
exponent at least $p$, and making us of the fact that if a variable is
required to have at least a certain degree, this is equivalent to simply
lowering the total degree of the monomial.

The asserted formula for $N^D_p(n,s)$ also follows easily from the 
definitions, since the count may be split up over
the choice of which $D$ variables have degree less than $p/2$, and then the
desired monomial is a product of a monomial with those variables, each with 
degree less than $p/2$, with a monomial of the remaining variables, each 
with degree greater than $p/2$. Summing over possible degrees of the two
separate monomials gives the formula.
\end{proof}

\begin{prop}\label{con-constant} We can classify completely all kernel map
classes 
with kernel isomorphic to $\sO(-m) \oplus \sO(m-n+\delta)$ which can be made 
via transport to have either $g_{11}$ or $g_{12}$ equal to $0$.
We may describe them as (note that despite the geometry in the
description, we make no claim of any {\it a priori} scheme or variety 
structure):
\begin{nlist}
\itm there are $N_p(n, mp+d)$ transport-antidiagonalizable classes.
\itm For each $D$, there are $N_p^D(n, mp+d)$ $\P^{D-2-n+\delta+2m}$'s of 
classes of non-transport-antidiagonalizable kernel maps for which $g_{11}$
may be transported to $0$.
\itm If $m \neq n-\delta-m$, there are an additional (distinct)
$N_p(n, (n-\delta-m)p+d)$ transport-diagonalizable classes.
\itm Again if $m \neq n-\delta-m$, for each $D$ there are an additional 
(distinct) $N_p^D(n, (n-\delta-m)p+d)$ $\P^{D-1}$'s of classes of 
non-transport-diagonalizable kernel maps for which $g_{12}$ may be 
transported to $0$.
\end{nlist}

In particular, in the case $mp+\delta< d$, all possible kernel maps are 
classified by (1) and (2).
\end{prop}

\begin{proof}
We begin with the case that $g_{11}=0$. Scaling as necessary, we have
$g_{21}g_{12} = \prod _{i=1} ^n (x-\lambda_i)^p$, so fix the orders at each
$\lambda_i$ of $g_{21}$ (equivalently, $g_{12}$). There are $N_p(n, mp+d)$
choices, by definition. Clearly, we get only one antidiagonalizable one
given the choices of orders. We next examine the
non-antidiagonalizable ones. Let $D$ be the number of $i$ such that
$g_{21}$ has order less than $p/2$ at $\lambda_i$; this will be the number
of $c_i$ which are uniquely determined under our criterion of Proposition
\ref{con-crit2}. Indeed, for such $c_i$, noting that $p-\ord_{\lambda_i}
g_{21} = \ord_{\lambda_i} g_{12}$, this criterion
tells us precisely that for each $\lambda_i$, there is a $c_i$ such that
$g_{22} \equiv c_i g_{21} \pmod{(x-\lambda_i)^{\ord_{\lambda_i}g_{12}}}$. 
In the cases that $c_i$ can be arbitrary, we may for convenience consider 
them to be $0$. We then observe that if we choose values for the remaining 
$c_i$, there is at most one transport-class with those values, 
since the Chinese remainder theorem says $g_{22}$ is is uniquely determined 
modulo $g_{12}$ by its values modulo
$(x-\lambda_i)^{\ord_{\lambda_i}g_{12}}$ for all $i$.
Now, there are $D$ of the $c_i$ which must be
specified, and they cannot all be the same, since in that case one could
arrange by a single column operation for $g_{12}$ to divide $g_{22}$, which
then means we are in the
transport-antidiagonalizable case. Moreover, by Proposition \ref{con-cimove}
we can do a global column operation to set the first $n-\delta-2m+1$ of 
the $c_i$ to $0$ (since powers of distinct numbers are always linearly
independent), reducing us to $D-n+\delta+2m-1$ choices, and we can also 
scale all the remaining $c_i$. So, we have a $\P^{D-2-n+\delta+2m}$
of distinct choices for the $c_i$, each corresponding to a unique class of
kernel maps. When $m<n-\delta-m$, from Proposition \ref{con-cimove} we 
know these are the only possibilities, so we are done. On the other hand, 
when $m=n-\delta-m$, 
we note that the only possibilities for right transport which preserve
$g_{11}=0$ are the upper triangular ones, which correspond precisely to the
translation and scaling we have already used, so this case works out exactly
the same way. This finishes cases (1) and (2).

For (3) and (4), first note that when $m=n-\delta-m$, one can globally 
switch columns, so the classes with $g_{12}$ transportable to $0$ are the
same as the ones we have already classified. For $m<n-\delta-m$, they are
distinct, since if either of the $g_{1j}$ is $0$, it is clear no 
transport-equivalent matrix could have the other $0$ instead. Thus, we argue 
in exactly the same way in this case, except that
for convenience we classify kernel map classes by the $c_i$ for the mirror
criterion, and we also have to note that globally in this case we cannot 
translate the $c_i$ at all, since any non-trivial column operation would make 
$g_{12} \neq 0$, so we get a $\P^{D-1}$ rather than a 
$\P^{D-2-n+\delta+2m}$. 
\end{proof}

\begin{rem} With this proposition, we already see polynomials in $p$ arising 
in counting connections with
a fixed set of poles on a fixed vector bundle. Ultimately, the numbers of
this proposition will not come into the calculation of the number of
connections we are interested in for the Frobenius-unstable vector bundles
of \cite{os11}, but that number
will also be a polynomial in $p$, strongly suggesting the existence of a
more general underlying phenomenon.
\end{rem}

\section{Maps from $\P^1$ to $\P^1$}\label{s-con-maps}

Continuing with the notation of the previous section, we have fully 
analyzed classes of kernel maps in which one of $g_{11}$ or
$g_{12}$ may be transported to $0$. To analyze the remaining kernel maps,
we shift focus considerably. We will associate a rational function $f_g$ to
each kernel map class, and examine the induced correspondence to complete
our general classification of logarithmic connections with vanishing
$p$-curvature, concluding in particular the statement of Theorem
\ref{con-main}.

\begin{warn}In order to streamline the proofs in this section, whenever we
refer to the $c_i$ or criterion of Proposition \ref{con-crit2}, we will mean 
the mirror criterion under which scalar multiples of the right column are 
added to the
left. 
\end{warn}

We begin with some notation and observations: first, since $\det(S)$ is
supported at the $\lambda_i$, the GCD of the coefficients of $S$ must
likewise be.

\begin{notn}\label{con-map-notn} Set $\alpha_i$ so that the GCD of the 
coefficients of $S$ is
$\prod(x-\lambda_i)^{\alpha_i}$. Factoring this out from $S$, write
$\hat{S}=(\hat{g}_{ij})$ for the resulting matrix, whose coefficients have 
no nontrivial common divisor. Now, let
$g_1$ be the GCD of $\hat{g}_{11}$ and $\hat{g}_{12}$, write $\beta_i =
\ord_{\lambda_i} g_1$, and finally write $f_g := \frac{g_{12}}{g_{11}},$
considered as an endomorphism of $\P^1$.
\end{notn}

We make the following observation: formally locally at each $\lambda_i$, we 
can transport-diagonalize $S$ to have powers of $t_i$ on the diagonal,
obtaining
two positive integers summing to $p$ as the exponents. Momentarily writing
$\alpha'_i$ for the lesser of the two, we note that $t_i^{\alpha'_i}$ is the
GCD of the coefficients of the diagonalized matrix, and since GCDs are 
unchanged by multiplication by invertible matrices, it must also have been 
the GCD of the coefficients of $S$ (over $k[[t_i]]$); hence, $\alpha'_i =
\alpha_i$. 

Note that since we assumed $\delta p -d < d$, the $g_{1j}$, and in particular, $f_g$, are unaffected by left transport. We also see easily that one of $g_{11}, g_{12}$ may be
transported to $0$ if and only if $S$ is transport-equivalent to a
kernel map with $f_g$ having degree $0$, hence constant. Thus:

\begin{cor}The kernel map classes classified in Proposition 
\ref{con-constant} are precisely those for which the associated 
endomorphism $f_g:\P^1
\rightarrow \P^1$ can be made constant under transport.
\end{cor}

We also note that via the constant row and column operations available to us
under global transport, we can without loss of generality assume we are in
the following situation. 

\begin{sit}\label{con-map-sit} We have normalized so that $\ord _{\lambda_i}
g_{22} = \alpha_i$ for all $i$, $\ord_{\lambda_i} g_{12} = \alpha_i + \beta_i$,
and $\deg g_{12} = (n-m)p-d$.
\end{sit}

We now analyze the situation further:

\begin{prop} We have $\beta_i \leq p-2\alpha_i$ for all $i$; moreover, $f_g$ has degree $(n-m)p-d-\sum_i \alpha_i - \sum_i \beta_i$,
and is ramified to order at least $p-2\alpha_i-\beta_i$ at each
$\lambda_i$, and $(n-\delta-2m)p$ (when this is non-zero) at infinity. 
\end{prop}

\begin{proof} The inequality $\beta_i \leq p-2\alpha_i$ is necessary for the determinant to have order $p$ at $\lambda_i$. Next, by definition 
$f_g$ has degree $(n-m)p-d-\sum_i \alpha_i
- \deg g_1$. Noting that $g_1$ will divide the determinant of $\hat{S}$, it
must be supported at the $\lambda_i$, so we also have $\deg g_1 = \sum_i
\beta_i$. Examining 
the (mirror) criterion of Proposition \ref{con-crit2}, we see that the
requirement
that $(g_{11}-c_i g_{12})(g_{22})$ vanish to order at least $p$ at 
$\lambda_i$, since we had arranged for $\ord_{\lambda_i} \hat{g}_{22} = 0$,
gives the desired ramification condition at $\lambda_i$. Finally, the 
ramification at
infinity follows because we had set $\deg g_{12} = (n-m)p-d$, so 
it has degree at least $(n-\delta-2m)p$ greater than $g_{11}$.
\end{proof}

In particular, we see that when $f_g$ is nonconstant and $\beta_i \neq p - 2\alpha_i$, the $c_i$ are 
all uniquely determined as $\frac{1}{f_g(\lambda_i)}$. In the case that 
$\beta_i = p - 2\alpha_i$, we need to additionally specify the $c_i$, but 
in applications (for instance, in Theorem \ref{con-main}), this situation 
will not arise. In determining necessary and 
sufficient conditions to fill in the $\hat{g}_{2j}$ from the 
$\hat{g}_{1j}$ in such a way as to satisfy our criterion, we find:

\begin{prop}\label{con-fillin} For any given choice of $\alpha_i$ and $\hat{g}_{1j}$ as
prescribed by the previous proposition, and for any choice of $c_{i_\ell}\in k$ for each $\beta_{i_\ell} = p -2\alpha_{i_\ell}$, and $c_{i_\ell} \neq \frac{1}{f_g(\lambda_{i_\ell})}$, there is a 
unique corresponding kernel map class if and only if for all $i$ such that
$0 < \beta_i < p -2\alpha_i$, $f_g$ has precisely the minimum required
ramification, i.e. $f_g$ ramifies to order precisely $p-2\alpha_i-\beta_i$ at
$\lambda_i$. Otherwise, there will be no corresponding kernel map.
\end{prop}

\begin{proof} We first show that if the conditions on $\hat{g}_{1j}$ are
satisfied, we get a unique corresponding kernel map class:
that is, given $\hat{g}_{11}$ and
$\hat{g}_{12}$, there is a unique way (up to transport) to fill in
$\hat{g}_{21}$ and $\hat{g}_{22}$ which satisfies the (mirror) criterion of
Proposition \ref{con-crit2}. This will follow from standard results
on generators of ideals over PIDs: we need to choose the bottom
row so that the determinant is $\Delta = \prod (x-\lambda_i)^{p- 2\alpha_i}$;
the solutions to
$\hat{g}_{11} h_2 - \hat{g}_{12} h_1 = \Delta$ are expressible for some
particular choice of $h_1, h_2$ as $h_1 + q
\frac{\hat{g}_{11}}{g_1}, h_2 + q \frac{\hat{g}_{12}}{g_1}$ as $q$ varies 
freely.  In particular, two ways of filling in the bottom row are transport 
equivalent if and only if their corresponding $q$'s differ by a multiple of
$g_1$, so we will need to check that the criterion determines 
$q$ precisely modulo $g_1$. We also observe that given $q$ modulo $g_1$, we
can always choose a representative polynomial for it so that the resulting
$g_{21}, g_{22}$ have the correct degrees: changing $q$ by a multiple of
$g_1$ corresponds to subtracting a multiple of the first row from the
second, which can always, for instance, force the degree of $g_{22}$ to be 
strictly smaller than $(n-m)p-d < (n-\delta-m)p+d$, without changing 
the determinant, and this forces $g_{21}$ to have degree exactly $mp+d$. 
Note also that some $h_1, h_2$ as above must exist because we assume
that $\beta_i \leq p - 2 \alpha_i$. 

Now, note that any $q$ yields a solution satisfying the order conditions
along the diagonal, antidiagonal, and the top row: indeed, our ramification 
condition gives order at least $p$ along the top row and 
diagonal (after column operation by $c_i$), and the determinant then
forces the antidiagonal to also have order at least $p$ at all
$\lambda_i$. In particular, $\ord_{\lambda_i} h_1 - c_i h_2 \geq p- 2
\alpha_i - \beta _i$, since we arranged for 
$\ord _{\lambda_i} \hat{g}_{12} = \beta_i$ at all $i$.
Next, we know that if we can fill in 
the bottom row so to satisfy our criterion, we can do it with
$\hat{g}_{22}$ non-vanishing at all $\lambda_i$, and conversely, if
$\hat{g}_{22}$ is non-vanishing at all $\lambda_i$, our criterion
requires precisely (in addition to the determinant being correct) that
$\ord_{\lambda_i} (g_{21}- c_i g_{22}) g_{22} \geq p$, or equivalently,
$\ord _{\lambda_i} (\hat{g}_{21} - c_i \hat{g}_{22}) \geq p - 2 \alpha_i$.
Plugging in our expressions for possibilities for $\hat{g}_{21}$ and
$\hat{g}_{22}$, we get 
$$\ord_{\lambda_i} (h_1- c_i h_2 + q
(\frac{\hat{g}_{11}}{g_1}-c_i \frac{\hat{g}_{12}}{g_1})) \geq p - 2
\alpha_i.$$ 
Now, we observed earlier that $h_1 - c_i h_2$ has order at least $p-2 
\alpha_i - \beta_i$. There are three cases to consider. If we have
$\beta_i=0$, there is nothing to check. If $0 < \beta_i < p-2\alpha_i$, the
latter term above has order precisely $p- 2 \alpha _i - \beta_i$ by
hypothesis, 
in which case $q$ will be determined uniquely modulo $(x-\lambda_i)^{\beta_i}$ by our order condition. Finally, if $\beta_i = p - 2\alpha_i$ and further $c_i \neq \frac{1}{f_g(\lambda_i)}$, the ramification condition is irrelevant and we have that the order of the second term is again $p- 2 \alpha _i - \beta_i$, as in the previous case. In this last situation, we also check by solving for $q$ that different choices of $c_i$ necessarily yield different choices of $q$ modulo $(x-\lambda_i)^{\beta_i}$. Combining these for all $i$ by the Chinese remainder theorem determines a unique $q$ 
modulo $g_1$, giving us our unique kernel map class corresponding to $f_g$, 
as desired.

Conversely, if $g_1$ has support at $\lambda_i$, then either $f_g$ has
to ramify to precisely the required order at $\lambda_i$, and no higher, or we must have $\beta_i=p-2\alpha_i$ with 
$c_i \neq \frac{1}{f_g(\lambda_i)}$:
Since $g_1$ is supported at $\lambda_i$, we have $\ord_{\lambda_i} g_{12} >
\ord _{\lambda_i} g_{22}$, so under our 
criterion, after column translation, since $\ord_{\lambda_i} g_{21} g_{22}
\geq p$, we obtain $\ord _{\lambda_i} g_{12} g_{21} > p$, and the 
determinant condition then implies
that $\ord _{\lambda_i} g_{11} g_{22} = p$. If $\beta_i < p -2\alpha_i$, we saw that the $c_i$ was uniquely determined as $\frac{1}{f_g(\lambda_i)}$, 
meaning that we cannot have any extra ramification at $\lambda_i$. On the other hand, if $\beta_i = p -2\alpha_i$, we cannot have $\ord _{\lambda_i} g_{11} g_{22} = p$ unless $c_i \neq \frac{1}{f_g(\lambda_i)}$, giving us our desired restrictions.
\end{proof}

We now obtain the following theorem.

\begin{thm}\label{con-main-nonconst} Fix $d,n,m, \delta$ with $\delta=0
\text{ or }1$, such that $\delta p -d \leq d$ and $m \leq n-\delta-m$,
together with points $P_1, \dots, P_n$ on $\P^1$. Then transport equivalence
classes of connections $\nabla$ on $\sO(\delta p -d) \oplus \sO(d)$ on $\P^1$
having trivial determinant, vanishing $p$-curvature, and logarithmic poles
at the $P_i$, with the kernel of $\nabla$ isomorphic to $\sO(-mp) \oplus
\sO((m-n+\delta)p)$, are classified as follows: 
\begin{nlist}
\itm Connections not inducing a connection on $\sO(d) \subset \sO(\delta p -d)
\oplus \sO(d)$, and having residues at $P_i$ with eigenvalues $(\alpha_i,
-\alpha_i)$ for $0<\alpha_i<\frac{p}{2}$, are classified by equivalence
classes of triples $(f, \{\beta_i\}_i, \{c_{i_\ell}\}_\ell)$, where the
$\beta_i \in \Z_{\geq 0}$ are bounded by $p-2\alpha_i$, there is a
$c_{i_\ell}\in k$ for each $i_\ell$ with $\beta_{i_\ell}=p-2\alpha_{i_\ell}$,
and $f$ is a separable rational function on $\P^1$ of degree 
$(n-m)p-d-\sum_i \alpha_i - \sum_i \beta_i$, ramified to order at least
$p-2\alpha_i - \beta_i$ at each $\lambda_i$ (with equality whenever
$0<\beta_i<p-2\alpha_i$), and further mapping infinity to infinity to order 
at least
$(n-\delta-2m)p$. Finally, we require $c_{i_\ell} f(\lambda_{i_\ell}) \neq 1$ for all $\ell$. The equivalence relation is generated by fractional linear transformation of $f$, and translation by inseparable polynomials 
$f_0$ of degree $\leq (n-\delta-2m)p$, with the $c_{i_\ell}$ related by the same linear fractional transformation or by $f_0(\lambda_{i_\ell})$, as appropriate.
\itm Connections inducing a connection on $\sO(d) \subset \sO(\delta p -d)
\oplus \sO(d)$ are classified in two categories. The first are those
classified by Proposition \ref{con-constant}. Connections in the second
category, having residues at $P_i$ with eigenvalues $(\alpha_i, -\alpha_i)$
for $0<\alpha_i<\frac{p}{2}$ satisfying $d+\sum_i \alpha_i < (m+\delta)p$, 
are classified by equivalence classes of pairs $(f, \{c_{i_\ell}\}_\ell)$, 
with the $i_\ell$ some subset of $\{1,\dots,n\}$, and
$f$ an inseparable rational function on $\P^1$, of degree $(n-m)p-d-\sum_i
\alpha_i-\sum_{\ell} (p-2\alpha_{i_\ell})$, and mapping infinity to infinity
to order at least $(n-\delta-2m)p$. As before, we require $c_{i_\ell}
f(\lambda_{i_\ell}) \neq 1$ for all $\ell$, and the equivalence relation is
the same as above. 
\end{nlist} 
\end{thm}

\begin{proof} We begin by noting that the
hypothesis that the connections in question do not induce a connection on
$\sO(d)$ is equivalent to $f_g$ being nonconstant and separable, 
since this is precisely when the upper right coefficient in Equation 
\ref{con-2x2T} will be non-zero. To see that this is equivalent to
restricting to separable $f_g$ in the ``non-constant'' case described by
Theorem \ref{con-main-nonconst}, it suffices to observe that if a kernel map
is transport equivalent to one with $f_g$ constant, then its $f_g$ is
necessarily inseparable. 

Now, it is easy to see that in the case $m=n-\delta-m$, we get each
kernel map class corresponding to a unique function, with transport
corresponding to automorphism of
$\P^1$. To see how the $c_{i_\ell}$ for $\beta_{i_\ell}=p-2\alpha_{i_\ell}$ change under such an automorphism, it suffices to note that although they are not determined by $f_g(\lambda_i)$, it follows from the proof of the previous proposition that they are determined as $\frac{g_{21}}{g_{22}}(\lambda_i)$, and thus change by the same automorphism. In the case $m<n-\delta-m$, it is clear from the definition of
$f_g$ that transport of a kernel map can change $f_g$ precisely by an 
inseparable polynomial of degree at most $(n-\delta-2m)p$. Thus, the condition that $d+\sum_i \alpha_i < (m + \delta)p$ insures that $f_g$ is not transport-equivalent to a constant function. The translation of the $c_{i_\ell}$ in this situation is given by Proposition \ref{con-cimove}. Putting all this
together with the previous propositions, we conclude the statement of the theorem.
\end{proof}

We further show:

\begin{prop} In case of Theorem \ref{con-main-nonconst}, if also $m=n-\delta-m$ and $\beta_i=0$ for all $i$, the classification holds over 
$k[\epsilon]$.
\end{prop}

\begin{proof} We begin by remarking that in our situation, we know that the
kernel of the deformed connection is a deformation of the kernel of the
original connection, by Proposition \ref{con-ker-def}. In the case
$m=n-\delta-m$, our kernel bundle has no non-trivial deformations, so a
deformation of a connection simply gives a deformation of the $\vp$ of the
kernel map, leaving the $\F$ intact.

Thus, we may represent our connection over $k[\epsilon]$ as a kernel 
map given by a matrix $(g_{ij} + \epsilon h_{ij})$, where we continue
with the notation of Notation \ref{con-map-notn} for the kernel map over
$k$ given by $(g_{ij})$, and assume the
$g_{ij}$ have been normalized as in Situation \ref{con-map-sit}. 
We know from Corollary
\ref{con-ct2} that our kernel matrix must still be formally locally
diagonalizable with the same eigenvalues over $k[\epsilon]$, so our
observation that 
$\alpha_i$ was alternatively described as the smaller eigenvalue of
the formally locally diagonalized kernel map gives us that
each of the $h_{ij}$ must also vanish to order at least $\alpha_i$ at
$\lambda_i$, and we set $\hat{h}_{ij}:=\frac{h_{ij}}{\prod_i (x-
\lambda_i)^{\alpha_i}}$. Because we have assumed $\beta_i=0$, it follows 
that $\frac{\hat{g}_{12}+\epsilon \hat{h}_{12}}{\hat{g}_{11}+\epsilon
\hat{h}_{11}}$ is a deformation of $f_g$ maintaining the same degree. It is
easy to check that the fact that Proposition \ref{con-crit2} holds over
$k[\epsilon]$ allows the same analysis as before to show that our
deformation preserves the required ramification, and it is clear that
transport still corresponds to postcomposition by an automorphism of $\P^1$. 

It therefore remains only to show that given an appropriate deformation of
$f_g$, we can still uniquely produce a corresponding kernel map over
$k[\epsilon]$. We therefore suppose we are given $\alpha_i$ for each $i$, 
$\hat{g}_{11}+\epsilon \hat{h}_{11}$, and $\hat{g}_{12}+\epsilon 
\hat{h}_{12}$. We may further suppose that we have $\hat{g}_{21}$ and
$\hat{g}_{22}$ satisfying the required determinant, degree, and vanishing
conditions modulo $\epsilon$, so we are simply trying to uniquely produce
$\hat{h}_{21}, \hat{h}_{22}$ to do likewise over $k[\epsilon]$. We first
consider the determinant condition: with $\hat{h}_{21}= \hat{h}_{22}=0$, the
determinant will be off by $\epsilon(\hat{g}_{22}\hat{h}_{11} - \hat{g}_{21}
\hat{h}_{12})$ from the desired $\prod _i (x - \lambda_i)^{p-2\alpha_i}$. We
see that we want to choose $\hat{h}_{21}, \hat{h}_{22}$ so that we have
$\hat{g}_{11}\hat{h}_{22} - \hat{g}_{12} \hat{h}_{21} =
\hat{g}_{22}\hat{h}_{11} - \hat{g}_{21} \hat{h}_{12}$, and this will be
possible if and only if $g_1 := \gcd(\hat{g}_{11}, \hat{g}_{12}) | 
(\hat{g}_{22}\hat{h}_{11} - \hat{g}_{21} \hat{h}_{12})$. However, since we
have assumed that all $\beta_i =0$, we have $g_1 = 1$, and may choose
$\hat{h}_{21}, \hat{h}_{22}$ to give the desired determinant. Moreover, 
given any fixed way of filling in the bottom row to give the right 
determinant, we see that all possible choices (with the same $\hat{g}_{21},
\hat{g}_{22}$) are given precisely as those obtained by adding 
$\epsilon$-multiples of the top row to the bottom, which gives the desired 
uniqueness. We can then use the same
argument as in the proof of Proposition \ref{con-fillin} to force the
degrees of $\hat{h}_{21}, \hat{h}_{22}$ to be bounded by $mp+d-\sum_i
\alpha_i, (n-\delta-m)p+d - \sum_i \alpha_i$ as required. Lastly, 
we must verify the vanishing condition imposed by Proposition 
\ref{con-crit2}; since everything will be multiplied through by 
$\prod_i (x - \lambda_i)^{\alpha_i}$, it is enough to verify that at
each $\lambda_i$, after column operation by $c_i$, we will have
$\min\{\ord_{\lambda_i} (\hat{g}_{11} + \epsilon \hat{h}_{11}),
\ord_{\lambda_i} (\hat{g}_{21} + \epsilon \hat{h}_{21})\} \geq p-2\alpha_i.$
By hypothesis, since $\ord _{\lambda_i} \hat{g}_{22} = 0$, this will already 
be satisfied modulo $\epsilon$, and the
ramification condition gives precisely $\ord_{\lambda_i} (\hat{h}_{11}) \geq 
p-2\alpha_i$, so it remains only to check that $\ord_{\lambda_i} (\hat{h}_{21})
\geq p-2\alpha_i$. However, we now see that  
$\ord_{\lambda_i} \hat{g}_{11}\hat{h}_{22} - \hat{g}_{12} \hat{h}_{21} =
\ord_{\lambda_i} \hat{g}_{22}\hat{h}_{11} - \hat{g}_{21} \hat{h}_{12} \geq 
p-2\alpha_i$, so since $\ord_{\lambda_i} \hat{g}_{11} \geq  p - 2\alpha_i$
and $\ord_{\lambda_i} \hat{g}_{12} =0$, we get the desired inequality.
\end{proof}

\begin{rem} The condition that $\beta_i=0$ in the above proposition is
unnecessary if one is willing to look at $g^1_d$'s rather than maps, 
and do slightly more analysis of vanishing conditions. However, we will
only need the case $\beta_i=0$. 
\end{rem}

We are now in a position to give:

\begin{proof}[Proof of Theorem \ref{con-main}] We note that the degree
and ramification conditions imposed in Theorem \ref{con-main-nonconst}, by
the separability of $f_g$ and the Riemann-Hurwitz formula, mean that no
additional ramification can occur and the $\beta_i$ must all be zero.
It therefore suffices to show that the only case which can actually occur
when the $P_i$ are general is the case $n-\delta-m=m$.

Now, suppose that $n-\delta-m > m$; since $n=2d+2$ is even, we must have 
$n-\delta-2m \geq
2$, and we see from Riemann-Hurwitz that there are two cases to consider:
either the ramification at infinity is exactly $(n-\delta-2m)p$, or it is
$(n-\delta-2m)p+1$. For the former case, if we subtract an
appropriate multiple of $x^{(n-\delta-2m)p}$, the Riemann-Hurwitz formula 
and the ramification at $\infty$ shows that we must
reduce the degree and ramification index at $\infty$ by precisely $1$. 
Noting that our choice of the point $\infty$ was arbitrary, both 
possibilities are then ruled out for general $P_i$ by \cite[Prop. 5.4]{os7}.
We conclude that for $P_i$ general, $n-\delta-m=m$ is
the only case that occurs, as desired. Finally, given this, the previous 
proposition shows that the classification still holds for first-order 
infinitesmal deformations.
\end{proof}

\begin{rem} We see in particular that the relationship between connections
and maps really is more complicated in the case of more than three
poles/ramification points, and one cannot hope to treat it as generally as
Mochizuki treated the three-point case; specifically, we see that
connections with $m \neq n-\delta-m$, which is to say those corresponding 
to maps 
with additional ramification at infinity will deform, as the poles move, to 
connections with $m = n-\delta-m$, which are lower-degree maps. This is 
a result of the fact that Grothendieck's splitting theorem for locally free 
sheaves on $\P^1$, used in an essential way in our argument, only holds 
over a field.
\end{rem} 

\bibliographystyle{hamsplain}
\bibliography{hgen}
\end{document}